\documentclass[10pt]{aiaa-tc}
\usepackage{epsfig,dcolumn,portland,amsmath,amssymb}
\usepackage{setspace}
\begin{document}
\bibliographystyle{plainnat} 

\title{Real-Time First Order Guidance Strategies for Trajectory Optimization in UAVs by Utilizing Wind Energy}
\author{Kamran Turkoglu 
\thanks{$^{1}$Kamran Turkoglu is with Faculty of Aerospace Engineering, College of Engineering,
        San Jose State University, One Washington Square, San Jose, CA 95192, USA
        {\tt\small kamran.turkoglu@sjsu.edu}}%
}

\maketitle

\begin{abstract}                          
This paper presents real-time guidance strategies for unmanned aerial vehicles (UAVs) that can be used to enhance their flight endurance by utilizing {\sl insitu} measurements of wind speeds and wind gradients. In these strategies, periodic adjustments would be made in the airspeed and/or heading angle command for the UAV to minimize a projected power requirement at some future time. In this paper, UAV flights are described by a three-dimensional dynamic point-mass. Onboard closed-loop trajectory tracking logics that follow airspeed vector commands are modeled using the method of feedback linearization. A generic wind field model is assumed that consists of a constant term plus terms that vary sinusoidally with respect to the location. To evaluate the benefits of these strategies in enhancing UAV flight endurance, a reference strategy is introduced in which the UAV would seek to follow the desired airspeed in a steady level flight under zero wind. A performance measure is defined as the average power consumption both over a specified time interval and over different initial heading angles of the UAV. A relative benefit criterion is then defined as the percentage improvement of the performance measure of a proposed strategy over that of the reference strategy. Extensive numerical simulations are conducted. Results demonstrate the benefits and trends of power savings of the proposed real-time guidance strategies.
\end{abstract}


\section{Introduction}

Because of the generally light weights and/or small sizes of unmanned aerial vehicles (UAVs), wind can play an important, sometimes crucial, role in their flight endurance and performance. Given the ubiquitous nature of the wind, it is highly desirable to devise UAV flight strategies that enable them to benefit from the wind.

However, a main challenge in utilizing wind energies for practical UAV flights is the need to obtain accurate and timely wind information in real time. Because UAVs may travel to remote regions where ground support systems are not readily available and typically regional wind field information is unknown, airborne measurements of the local wind information are essential. Therefore, an interesting problem of UAV wind utilization is to develop onboard guidance and control strategies that can take advantage of wind energies based on airborne measurements (or estimates) of the local winds.

This paper presents such real-time guidance strategies in which optimal adjustments are made to the airspeed and heading angle commands to minimize a projected power consumption at some future time to prolong a UAV flight, based on the current local wind conditions. The onboard feedback control system then tracks these modified commands. This process is repeated periodically.

There are pioneering works in the area of UAV flights utilizing wind energies. The developments and flight tests of practical guidance strategies for detecting and utilizing thermals by Allen\cite{allen:2005,allen:2007} and Edwards\cite{Edwards:2008} have illustrated the feasibility of these concepts. Boslough\cite{Boslough:2002} demonstrated the benefits of utilizing wind gradients through dynamic soaring using radio-controlled UAVs. Patel and Kroo\cite{patel:2006} studied the effect of wind in determining optimal flight control conditions under the influence of atmospheric turbulence. Langelaan and Bramesfeld\cite{langelaan:2008} studied how to exploit energy from high frequency gusts in the vertical plane for UAVs. Wharington\cite{wharington:1998,wharington:2004} presented methods for learning the wind patterns, based on local sensing and an appropriately selected reward function, and to fly most efficiently. Pisano\cite{Pisano:2009} investigated the gust sensitivity on the UAV dynamics as a function of aircraft size. In addition, Chakrabarty and Langelaan\cite{Langelaan:2009} presented a method for minimum energy path planning in complex wind fields using a predetermined energy map.  Lawrence and Sukkarieh\cite{Sukkarieh:2009} developed a framework for an energy-based path planning that utilizes local wind estimations for dynamic soaring. Rysdyk\cite{Rysdyk:2007} studied the problem of course and heading changes in significant wind conditions. McNeely\cite{McNeely:07} and et al. studied the tour planning problem for UAVs under wind conditions. McGee and Hedrick\cite{McGee:2007} presented a study of optimal path planning using a kinematic aircraft model.

Dynamic optimization methods have also been used to determine the full potential benefits of wind energy utilization when a regional wind model is known. Sachs, Knoll \& Lesch\cite{Sachs:1991} et al. studied optimal glider dynamic soaring in a wind gradient. Zhao and Qi\cite{zhao:2004b,qi:2005,zhao:2004} showed that under appropriate conditions with a full knowledge of the wind field in a region, a UAV can greatly enhance its endurance by properly utilizing wind energies. In addition to favorable wind patterns, such as wind gradient and thermals, Zhao\cite{zhao:2009} recently showed that even downdraft wind can be utilized to improve UAV performances. These results indicate that the utilization of wind energies for enhancing UAV flights is highly promising. Furthermore, dynamic optimization methods have also been applied to glider flights in winds\cite{Pierson:1980, Pierson:1978, Lorenz:1985}. In addition, Mueller, Zhao \& Garrard\cite{MuellerZhao:2008} studied optimal airship ascent flights by utilizing wind energy. In these studies, the nonlinear dynamic optimization formulation typically requires knowledge of regional winds and an iterative solution process, and thus may not be feasible for generating real-time guidance strategies. Still, they are useful in understanding fundamental patterns of optimal UAV flights in winds and providing benchmark results that can be used to evaluate real-time wind utilization strategies. Approximate solutions may also be obtained for nonlinear dynamic optimization problems to derive real-time guidance strategies.

Compared with the dynamic optimization studies, the current paper presents real-time guidance strategies that use {\sl insitu} wind measurements alone with no regional wind information, to reduce power consumptions. These strategies periodically adjust airspeed vector commands to take advantage of changes in the mean wind profile. Wind energies in the changing mean wind profile are generally of lower frequency compared with gust energies. As a result, the proposed real-time guidance strategies complement the previous works on real-time guidance and control methods that utilize gust energies.

In the current paper, only airspeed and heading adjustments are examined with a zero flight path angle command in order to evaluate the benefits of unconstrained guidance strategies. These strategies can be easily expanded to adjust the flight path angle as well, but this requires the incorporation of altitude boundary control in order to prevent the UAV from hitting the ground. Guidance strategies that need to respect flight constraints in altitude as well as over the horizontal region shall be reported later.

In the rest of the paper, three-dimensional point-mass equations are used to describe UAV motions in winds. Optimal adjustments for airspeed and heading angle are derived by minimizing a power consumption projected into the future. Models of closed-loop trajectory tracking are developed to follow airspeed vector commands, using the technique of feedback linearization. A generic wind pattern consisting of a constant term plus spatially varying terms is used to evaluate the average power consumption of the proposed strategies. In order to eliminate the impact of wind directions on the relative benefits of proposed strategies in UAV flights, the power consumption is averaged both over a specified time interval and over different initial heading angles of the UAV. This average is then compared with that of a reference strategy, in which the UAV seeks to maintain a steady level flight with the airspeed that would maximize the endurance in zero wind. Conclusions are drawn at the end.

\section{Equations of Motion and Constraints}

For the purposes of developing guidance strategies, UAV flights are represented by a dynamic point-mass model. The corresponding normalized equations of motion for a propeller-driven UAV are listed below\cite{jackson:1999}, where the UAV mass is assumed to be constant.
\subsection{Normalized Equations of Motion}

In order to increase numerical efficiency in the simulation studies, the above equations of motion are normalized by specifying a characteristic speed $V_n$ and mass $m$. We have
\begin{equation}
\begin{split}
\bar{V} = {V \over V_n},(\bar{x},\bar y,\bar{h}) = {({x,y, h}) \over V^2_n/g},\bar P = {P \over mg V_n},\bar{t} = {t \over V_n/g}
\end{split}
\end{equation}
\begin{equation}
( \ )' = {d ( \ ) \over d \bar{t}} = {V_n \over g} {d ( \ ) \over d \bar{t}}
\end{equation}
and
\begin{equation}
\bar{\rho} = { \rho S V^2_n\over 2 mg} = {\rho V^2_n  \over 2 (mg/S) }
\end{equation}
where the normalized air-density, $\bar{\rho}$, represents the combined effect of air density ($\rho$) and wing loading ($mg/S$) on UAV flights. Specifically, a smaller $\bar{\rho}$ corresponds to a larger wing loading (a heavier UAV) and/or thinner air, where a larger $\bar{\rho}$ represents a lighter UAV (with a smaller wing loading) and/or thicker air.

Using these normalizations, the normalized drag and lift become
\begin{equation}
\bar{D} = {D \over mg} = \bar{\rho} \bar{V}^2 C_D \hskip 10 pt
\bar{L} = {L \over mg}  = \bar{\rho} \bar{V}^2 C_L
\end{equation}
The normalized wind components are defined as
\begin{equation}
\bar W_{x,y,h} = {W_{x,y,h}\over V_n},~~~
         \bar W_{V,\gamma, \Psi} '  = {\dot W_{V,\gamma,\Psi}\over g},~~~
         \bar W_{x,y,h}'  = {\dot W_{x,y,h} \over g }
\end{equation}
where the normalized rates of wind speeds follow similar expressions as in Eqs. (7)-(9).

Then the set of normalized equations of motion are obtained as follows, where functional dependences of the wind terms are shown for convenience
\begin{equation}\label{VdotNormalized}
\begin{split}
&\bar{V}'  = {\bar{P} \over \bar{V}} - \bar{\rho} \bar{V}^2 (C_{D_0} + K C^2_L) - \sin\gamma - \bar W_V'(\bar{V}, \Psi, \gamma, \bar{x}, \bar{y}, \bar{h})
\end{split}
\end{equation}
\begin{equation}\label{PsidotNormalized}
\begin{split}
& \Psi' = {{\bar{\rho} \bar V C_L}\over\cos\gamma}\sin\mu - {1\over \bar V \cos\gamma} \bar W_\Psi'(\bar{V}, \Psi, \gamma, \bar{x}, \bar{y}, \bar{h})
\end{split}
\end{equation}
\begin{equation}\label{GammadotNormalized}
\begin{split}
& \gamma'  = \bar{\rho} \bar{V} C_L\cos\mu- {\cos\gamma\over \bar{V}}+ {1 \over{\bar{V}}} \bar W_\gamma'(\bar{V}, \Psi, \gamma, \bar{x}, \bar{y}, \bar{h})
\end{split}
\end{equation}
\begin{equation}\label{XdotNormalized}
\bar{x}'  = \bar{V} \cos\gamma \sin\Psi + \bar{W}_x(\bar{x}, \bar{y}, \bar{h})
\end{equation}
\begin{equation}\label{YdotNormalized}
\bar{y}'  = \bar{V} \cos\gamma \cos\Psi + \bar{W}_y(\bar{x}, \bar{y}, \bar{h})
\end{equation}
\begin{equation}\label{HdotNormalized}
\bar{h}'  = \bar{V} \sin\gamma + \bar{W}_h(\bar{x}, \bar{y}, \bar{h})
\end{equation}
Constraints on states and controls can also be expressed using normalized values.

\section{Problem Statement}

The main aim of this paper is to develop real-time guidance strategies to enhance the endurance of UAV flights based on {\sl insitu} wind information. Ideally, if the regional wind information is completely known in advance, optimal flight planning can be used to determine UAV flight trajectories that minimize the total power consumption over a specified time interval, subject to various constraints. However in this paper, it is assumed that only wind information at the current location of the UAV at the current time is available. This information includes values of wind speeds as well as wind gradients.

In general, different guidance strategies may be grouped into three basic categories:  {\bf action strategy}, {\bf velocity strategy}, and {\bf trajectory strategy}. This paper studies {\bf velocity guidance strategies} that can utilize {\sl insitu} wind information to enhance UAV endurance. A propeller-driven UAV is assumed, for which maximum endurance corresponds to minimum power consumption\cite{anderson:89}. Therefore, we seek to determine incremental adjustments in airspeed and heading angle to minimize the power consumption projected sometime into the future. Mathematically,
\begin{equation}\label{ProbDef}
\min_{\Delta \bar{V}_c, \Delta \Psi_c} \; \bar{I} = \bar{P}(\bar{t}_0+\Delta \bar{t})
\end{equation}
subject to all applicable constraints. Then, the UAV will be directed to track $\bar{V}_0 + \Delta \bar{V}_c$ and $\Psi_0 + \Delta \Psi_c$ commands.

Once adjustments in airspeed and heading angle are obtained, it takes some finite time for the UAV to achieve the desired changes via closed-loop tracking. As a result, a projected power consumption at $(\bar{t}_0+\Delta \bar{t})$ is used instead of the current power in Eq. (\ref{ProbDef}).

\section{Solution Strategies}

A key to solving the above problem is to develop an expression for the projected power required at the time $\bar{t}_0 +\Delta \bar{t}$, based on values of the current trajectory state variables as well as the current wind information. Therefore, $\Delta \bar{t}$ should neither be too large or too small. It needs to be larger than the typical settling time of inner closed-loop controls in order to ensure that any adjustments in airspeed and heading commands will have been achieved. At the same time, too large a $\Delta \bar{t}$ would reduce the accuracy of power consumption projections using the current state and wind information.

In this paper, it is assumed that the UAV intends to maintain a level flight: $L =W$, $\gamma = 0$, and $\bar{h}' = 0$.
\begin{equation}\label{CLinlevel}
C_L = {1 - \bar{W}'_\gamma \over \bar{\rho}\bar{V}^2 \cos\mu} \approx {1  \over \bar{\rho}\bar{V}^2 \cos\mu}
\end{equation}

From Eq. (\ref{VdotNormalized}), we have
\begin{equation}\label{power}
\begin{split}
\bar{P} & =\bar{V} \bar{V}'  +  \bar{\rho} \bar{V}^3 (C_{D_0} + KC_L^2)  + \bar{V} \bar{W}_V' \cr
	& \approx \bar{V} \bar{V}'  + \bar{\rho} \bar{V}^3 C_{D_0} + {K  \over \bar{\rho} \cos^2\mu} {1 \over \bar{V}} + \bar{V} \bar{W}_V'  \\
\end{split}
\end{equation}

Therefore, the projected power level at $\bar{t}_0+\Delta \bar{t}$ is given by
\begin{equation}\label{minPbar2}
 \bar{P}(\bar{t}_0+\Delta \bar{t})= \left[ \bar{V} \bar{V}'  + \bar{\rho} \bar{V}^3 C_{D_0} + {K  \over \bar{\rho} \cos^2\mu} {1 \over \bar{V}} + \bar{V} \bar{W}_V' \right]_{\bar{t}_0+\Delta \bar{t}}
\end{equation}

It is assumed that by the time $\bar{t}_0+\Delta \bar{t}$, any commanded changes in airspeed and heading angle will have been mostly achieved via closed-loop tracking. Therefore, the vehicle is basically in a steady state: $\bar{V}' \approx 0$ and $\mu \approx 0$. We have,

\begin{equation}\label{minPbar3}
\begin{split}
\bar{P}(\bar{t}_0+\Delta \bar{t}) \approx &~ \bar{\rho} [\bar{V}_0 + \Delta \bar{V}_c]^3 C_{D_0} + {K  \over \bar{\rho} } {1 \over \bar{V}_0 + \Delta V_c } \\
&~ + [\bar{V}_0 + \Delta \bar{V}_c] \bar{W}_V'(t_0+\Delta t)
\end{split}
\end{equation}
We now need to develop an expression for the $\bar{W}_V'(\bar{t}_0+\Delta \bar{t})$ term.

\begin{figure*}[!t]
\normalsize
\begin{equation}\label{Wvdot_looong}
\begin{split}
\bar{W}_V'= & \left( { \partial{\bar{W}_x} \over \partial{\bar{y}}} + { \partial{\bar{W}_y} \over \partial{\bar{x}}}\right) \bar{V} \cos^2\gamma \sin\Psi \cos\Psi
 + \left( { \partial{\bar{W}_x} \over \partial{\bar{h}}} + { \partial{\bar{W}_h} \over \partial{\bar{x}}}\right) \bar{V} \sin\gamma \cos\gamma \sin\Psi  + \left( { \partial{\bar{W}_y} \over \partial{\bar{h}}} + { \partial{\bar{W}_h} \over \partial{\bar{y}}}\right) \bar{V} \sin\gamma \cos\gamma \cos\Psi \\
&  + \left( { \partial{\bar{W}_x} \over \partial{\bar{x}}} \bar{V} \cos^2\gamma \sin^2\Psi \right) + \left( { \partial{\bar{W}_y} \over \partial{\bar{y}}} \bar{V} \cos^2\gamma \cos^2\Psi \right) + \left( { \partial{\bar{W}_h} \over \partial{\bar{h}}} \bar{V} \sin^2\gamma \right) + \left( { \partial{\bar{W}_x} \over \partial{\bar{\Delta}}} \cos\gamma \sin\Psi \right) + \left( { \partial{\bar{W}_y} \over \partial{\bar{\Delta}}} \cos\gamma \cos\Psi \right)  \\
& + \left( { \partial{\bar{W}_h} \over \partial{\bar{\Delta}}} \sin\gamma \right) \\
\end{split}
\end{equation}
\hrulefill
\vspace*{4pt}
\end{figure*}

\subsection{An Expression For Projected Wind Rate}
Because
\begin{equation}\label{windrate}
\bar{W}'_{( \ )} = {\partial \bar{W}_{( \ )} \over \partial \bar{x}} \bar{x}' + {\partial \bar{W}_{( \ )} \over \partial \bar{y}} \bar{y}'
+ {\partial \bar{W}_{( \ )} \over \partial \bar{h}} \bar{h}' + {\partial \bar{W}_{( \ )} \over \partial \bar{t}}
\end{equation}
where ${( \ )} = x, y, h$. By substituting previous expression in Eq.(\ref{windrate}), we have the expression given in Eq. (\ref{Wvdot_looong}), where
\begin{equation}
{(~~) \over \partial{\bar{\Delta}}} = \bar{W}_x { (~~) \over \partial{\bar{x}}} + \bar{W}_y { (~~) \over \partial{\bar{y}}} + \bar{W}_h { (~~) \over \partial{\bar{h}}} + { (~~) \over \partial{\bar{t}}}
\end{equation}

In level flights with negligible vertical winds, $\gamma = 0$, $\bar{W}_h \approx 0$,
\begin{equation}\label{WvBar2}
\begin{split}
\bar{W}_V' = & \left( { \partial{\bar{W}_x} \over \partial{\bar{y}}} + { \partial{\bar{W}_y} \over \partial{\bar{x}}}\right) \bar{V} \sin\Psi \cos\Psi \\
& + { \partial{\bar{W}_x} \over \partial{\bar{x}}} \bar{V} \sin^2\Psi + { \partial{\bar{W}_y} \over \partial{\bar{y}}} \bar{V}  \cos^2\Psi  \\
& + \left( \bar{W}_x { \partial{\bar{W}_x} \over \partial \bar{x} } + \bar{W}_y { \partial{\bar{W}_x} \over \partial \bar{y} } + { \partial{\bar{W}_x} \over \partial \bar{t} } \right) \sin\Psi \\
& + \left( \bar{W}_x { \partial{\bar{W}_y} \over \partial \bar{x} } + \bar{W}_y { \partial{\bar{W}_y} \over \partial \bar{y} } + { \partial{\bar{W}_y} \over \partial \bar{t} } \right) \cos\Psi \\
\end{split}
\end{equation}

Because only {\sl insitu} wind information is available, it is assumed that the current wind gradients shall stay constant over the immediate neighborhood around the current position of the UAV in the near future. This assumption shall be called the ``constant wind gradient assumption''.
Therefore, we obtain the final expression given in Eq.(\ref{WvBar2dt}).
\begin{figure*}[!t]
\normalsize
\begin{equation}\label{WvBar2dt}
\begin{split}
\bar{W}_V'(\bar{t}_0+\Delta \bar{t}) = & \left( { \partial{\bar{W}_x} \over \partial{\bar{y}}} + { \partial{\bar{W}_y} \over \partial{\bar{x}}}\right) (\bar{V}_0 + \Delta V_c)
\sin(\Psi_0 + \Delta \Psi_c) \cos(\Psi_0 + \Delta \Psi_c) + { \partial{\bar{W}_x} \over \partial{\bar{x}}} (\bar{V}_0 + \Delta V_c) \sin^2(\Psi_0 + \Delta \Psi_c) \\
& + { \partial{\bar{W}_y} \over \partial{\bar{y}}} (\bar{V}_0 + \Delta V_c) \cos^2(\Psi_0 + \Delta \Psi_c) + \bar{W}_x(\bar{t}_0+ \Delta \bar{t}) \left[ { \partial{\bar{W}_x} \over \partial \bar{x} } \sin(\Psi_0 + \Delta \Psi_c) + { \partial{\bar{W}_y} \over \partial \bar{x} } \cos(\Psi_0 + \Delta \Psi_c) \right] \\
& + \bar{W}_y(\bar{t}_0+ \Delta \bar{t}) \left[ { \partial{\bar{W}_x} \over \partial \bar{y} } \sin(\Psi_0 + \Delta \Psi_c) + { \partial{\bar{W}_y} \over \partial \bar{y} } \cos(\Psi_0 + \Delta \Psi_c) \right] + { \partial{\bar{W}_x} \over \partial \bar{t} } \sin(\Psi_0 + \Delta \Psi_c) \\
& + { \partial{\bar{W}_y} \over \partial \bar{t} }  \cos(\Psi_0 + \Delta \Psi_c) \\
\end{split}
\end{equation}
\hrulefill
\vspace*{4pt}
\end{figure*}
With the constant wind gradient assumption over $[\bar{t}_0, \bar{t}_0 + \Delta \bar{t}]$, we also have
\begin{equation}\label{windspeeds}
\begin{split}
\bar{W}_x (\bar{t}_0+\Delta \bar{t}) & \approx \bar{W}_{x_0} + { \partial{\bar{W}_x} \over \partial \bar{x} } \Delta \bar{x} + { \partial{\bar{W}_x} \over \partial \bar{y} } \Delta \bar{y} + { \partial{\bar{W}_x} \over \partial \bar{t} } \Delta \bar{t} \\
\bar{W}_y (\bar{t}_0+\Delta \bar{t}) & \approx \bar{W}_{y_0} + { \partial{\bar{W}_y} \over \partial \bar{x} } \Delta \bar{x} + { \partial{\bar{W}_y} \over \partial \bar{y} } \Delta \bar{y} + { \partial{\bar{W}_y} \over \partial \bar{t} } \Delta \bar{t} \\
\end{split}
\end{equation}

These expressions depend on $(\Delta \bar{x}, \Delta \bar{y})$, which depend on $\Delta \bar{V}_c$, $\Delta \Psi_c$, and $\Delta \bar{t}$ , and reciprocally on the wind components over the interval. Therefore, we need to develop expressions of $(\Delta \bar{x}, \Delta \bar{y})$ in order to complete the derivation of the projected wind rate expression.

\subsection{Expressions for Position Changes}
We now seek to develop expressions of $(\Delta \bar{x}, \Delta \bar{y})$ that show their dependencies explicitly on the increments of airspeed and heading angle. After experimenting with different methods, the following expressions are obtained. From Eqs. (\ref{XdotNormalized}) and (\ref{YdotNormalized}), we have for $\gamma=0$,
\begin{equation}\label{DxDyTrapezoidRule}
\begin{split}
\Delta \bar{x} = \bar{x}(\bar{t}_0+\Delta \bar{t}) - \bar{x}(\bar{t}_0) &  = \int_{\bar{t}_0}^{\bar{t}_0+\Delta \bar{t}} \left(\bar{V} \sin\Psi + \bar{W}_{x} \right) d \bar{t} \\
\Delta \bar{y} = \bar{y}(\bar{t}_0+\Delta \bar{t}) - \bar{y}(\bar{t}_0) &  = \int_{\bar{t}_0}^{\bar{t}_0+\Delta \bar{t}} \left( V \cos\Psi + \bar{W}_y \right) d \bar{t} \\
\end{split}
\end{equation}

Applying the trapezoidal rule\cite{Carnahan:1969} to integrate the above equations, with the assumption that both airspeed and heading angle will have achieved their commanded values at the end of the interval and the wind speeds are given by Eq. (\ref{windspeeds}), we obtain
\begin{equation}\label{DxDy}
\begin{split}
 {2 \over \Delta t } \Delta \bar{x} = &~ \bar{V}_0 \sin\Psi_0 + \bar{W}_{x_0} + (\bar{V}_0 + \Delta \bar{V}_c) \sin(\Psi_0 + \Delta \Psi_c) \\
 	& \quad + \bar{W}_{x_0} + { \partial{\bar{W}_x} \over \partial \bar{x} } \Delta \bar{x} + { \partial{\bar{W}_x} \over \partial \bar{y} } \Delta \bar{y} +
{ \partial{\bar{W}_x} \over \partial \bar{t} } \Delta \bar{t} \\
 {2 \over \Delta t } \Delta \bar{y} = &~ \bar{V}_0 \cos\Psi_0 + \bar{W}_{y_0} + (\bar{V}_0 + \Delta \bar{V}_c) \cos(\Psi_0 + \Delta \Psi_c) \\
 	& \quad + \bar{W}_{y_0} + { \partial{\bar{W}_y} \over \partial \bar{x} } \Delta \bar{x} + { \partial{\bar{W}_y} \over \partial \bar{y} } \Delta \bar{y} +
{ \partial{\bar{W}_y} \over \partial \bar{t} } \Delta \bar{t} \\
\end{split}
\end{equation}
Define,
\begin{equation}\label{B1B2}
\begin{split}
B_1 =&~ \bar{V}_0 \sin\Psi_0  + (\bar{V}_0 + \Delta \bar{V}_c) \sin(\Psi_0 + \Delta \Psi_c) \\
& + 2\bar{W}_{x_0} + { \partial{\bar{W}_x} \over \partial \bar{t} } \Delta \bar{t} \\
B_2 =&~ \bar{V}_0 \cos\Psi_0 + (\bar{V}_0 + \Delta \bar{V}_c) \cos(\Psi_0 + \Delta \Psi_c) \\
& + 2\bar{W}_{y_0} + { \partial{\bar{W}_y} \over \partial \bar{t} } \Delta \bar{t} \\
\end{split}
\end{equation}
we have
\begin{equation}
\left[ \begin{array}{cc} { 2 \over \Delta \bar{t}} - {\partial \bar{W}_x  \over \partial \bar{x} } & -{\partial \bar{W}_x \over \partial \bar{y} } \\
-{\partial \bar{W}_y \over \partial \bar{x} }  & { 2 \over \Delta \bar{t}} - {\partial \bar{W}_y \over \partial \bar{y} }  \\ \end{array} \right]
\left[ \begin{array}{c}  \Delta \bar{x} \\ \Delta \bar{y} \end{array} \right]
= \left[ \begin{array}{c} B_1 \\ B_2 \end{array} \right]
\end{equation}
or
\begin{equation}\label{DxDyComplete}
\begin{split}
\Delta \bar{x} & = {1 \over Q_{2D}} \left[ \left({ 2 \over \Delta \bar{t}}-{\partial{\bar{W}_{y}} \over \partial{\bar{y}}} \right) B_1 +  \left( {\partial{\bar{W}_{x}} \over \partial{\bar{y}}} \right) B_2 \right] \\
 \Delta \bar{y} & = {1 \over Q_{2D}} \left[ \left({\partial{\bar{W}_{y}} \over \partial{\bar{x}}} \right) B_1 +  \left( { 2 \over \Delta \bar{t}}-{\partial{\bar{W}_{x}} \over \partial{\bar{x}}}  \right) B_2 \right]
\end{split}
\end{equation}
where
\begin{equation}
\begin{split}
Q_{2D} =&~  \left( {2 \over \Delta \bar{t}} \right)^2 - {2 \over \Delta \bar{t}} \left( {\partial{\bar{W}_{x}} \over \partial{\bar{x}}}  +{\partial{\bar{W}_{y}} \over \partial{\bar{y}}} \right) \\
& + \left( {\partial{\bar{W}_{x}} \over \partial{\bar{x}}} {\partial{\bar{W}_{y}} \over \partial{\bar{y}}} - {\partial{\bar{W}_{x}} \over \partial{\bar{y}}} {\partial{\bar{W}_{y}} \over \partial{\bar{x}}} \right)
\end{split}
\end{equation}
For a sufficiently small update time-$\Delta t$, this expression shall always be nonzero; ensuring the existence of solutions for the position change expressions.

\subsection{Guidance Algorithms}

Based on the above derivations, we can now express the projected power consumption at $\bar{t}_0 + \Delta \bar{t}$ as a function of the current command adjustments in airspeed and heading angle. Then, the problem of reducing future power consumptions is to determine $\Delta \bar{V}_c$ and $\Delta \Psi_c$ from Eq.(\ref{StaticOptim}) with corresponding constraints in Eq.(\ref{StaticOptimConstr}).

\begin{figure*}[!t]
\normalsize
\begin{equation}\label{StaticOptim}
\min_{\Delta \bar{V}_c, \Delta \Psi_c} \bar{I}  = \bar{P}(\bar{t}_0 + \Delta \bar{t}) = \bar{I} (\Delta \bar{V}_c, \Delta \Psi_c; \Delta \bar{t}, \bar{X}_0)
\end{equation}
subject to
\begin{equation}\label{StaticOptimConstr}
\begin{split}
\max\{-\Delta \bar{V}_{c,\max}, \bar{V}_{\min} - \bar{V}_0\}  \leq \Delta \bar{V}_c \leq  \min \{ \Delta \bar{V}_{c,\max},  \bar{V}_{\max} - \bar{V}_0 \}  \\
\max\{-\Delta \Psi_{c,\max}, \Psi_{\min} - \Psi_0 \}  \leq \Delta  \Psi_c \leq \min \{ \Delta \Psi_{c,\max},  \Psi_c \leq \Psi_{\max} - \Psi_0 \} \\
\end{split}
\end{equation}
\hrulefill
\vspace*{4pt}
\end{figure*}

where $\Delta \bar{V}_{c,\max}$ and $\Delta \bar{\Psi}_{c, \max}$ are the maximum allowed incremental changes, and the initial state conditions required $\bar{X}_0$ include
\begin{equation}
\bar{X}_0 = \left\{ \bar{V}_0, \Psi_0, \bar{W}_{x_0}, \bar{W}_{y_0}, {\partial{\bar{W}_{x}} \over \partial{\bar{x}}}, {\partial{\bar{W}_{x}} \over \partial{\bar{y}}}, {\partial{\bar{W}_{y}} \over \partial{\bar{x}}}, {\partial{\bar{W}_{y}} \over \partial{\bar{y}}} \right\}
\end{equation}

As a reference strategy, it is assumed that the UAV follows a constant airspeed straight level flight. The airspeed is optimal in zero wind. In this case, the projected power expression in Eq. (\ref{minPbar3}) becomes
\begin{equation}
\bar{P} = \bar{\rho} \bar{V}^3 C_{D_0} + {K  \over \bar{\rho} } {1 \over \bar{V}}
\end{equation}
\begin{equation}
{ \partial \bar{P} \over \partial \bar{V}} = 0 \hskip 10pt \Rightarrow \hskip 10pt  \bar{V}^* = \left({ K \over 3  \bar{\rho}^2 C_{D_0}} \right)^{1 \over 4}
\end{equation}
This airspeed corresponds to the maximum endurance under zero wind in a steady level flight.

\section{First Order Adjustment Strategies}

Different algorithms can be used to solve the static optimization problem in Eq. (\ref{StaticOptim}). In this paper, first-order gradient algorithms are used to obtain solutions. In deriving a first-order gradient method, we approximate the projected power expression as
\begin{equation}\label{PowerApprox}
\bar{P}(\bar{t}_0+\Delta \bar{t}) \approx I_0 + \left( {\partial I \over \partial \Delta V_c}\right)_0  \Delta \bar{V}_c
+ \left( {\partial I \over \partial \Delta \Psi_c}\right)_0  \Delta \Psi_c
\end{equation}
where $( \ )_0$ corresponds to zero commanded adjustments or $(\Delta V_c = 0, \Delta \Psi_c = 0)$.

\subsection{Airspeed Adjustment Strategy}

A first-order incremental airspeed adjustment strategy can be determined from Eq. (\ref{PowerApprox}) as
\begin{equation}
\Delta \bar{V}_c = \begin{cases} 0 & \text{$|I_{\Delta \bar{V}_c}| < \epsilon$} \\
\eta \max\{-\Delta \bar{V}_{c,\max}, \bar{V}_{\min} - \bar{V}_0\}  &\text{$I_{\Delta \bar{V}_c} \ge \epsilon > 0$} \\
\eta \min \{ \Delta \bar{V}_{c,\max},  \bar{V}_{\max} - \bar{V}_0 \}  &\text{$I_{\Delta \bar{V}_c} \leq -\epsilon $}
\end{cases}
\end{equation}
where $\epsilon > 0$ is introduced to avoid numerical difficulties in implementation, and $\eta \in (0,1)$ is the adjustment stepsize. Expanding Eq. (\ref{minPbar3}) using Eq. (\ref{WvBar2}) leads to
\begin{equation}
\begin{split}
\bar{I}_{\Delta \bar{V}_c} & = 3 \bar{\rho} C_{D_0} \bar{V}^2_0 - {K \over \bar{\rho} \bar{V}^2_0} \\
& + \bar{V}_0 \sin(2\Psi_0) \left( {\partial \bar{W}_x \over \partial \bar{y}} + {\partial \bar{W}_y \over \partial \bar{x}}
+ {\partial \bar{W}_x \over \partial \bar{x}} \tan\Psi_0 + {\partial \bar{W}_y \over \partial \bar{y}} \cot\Psi_0 \right) \\
& + \left( {\partial \bar{W}_x \over \partial \bar{x} } \sin\Psi_0 + {\partial \bar{W}_y \over \partial \bar{x} } \cos\Psi_0 \right)
    \left[ \left( \bar{W}_x \right)_0 + \bar{V}_0 \left( {\partial \bar{W}_x \over \partial \Delta \bar{V}_c}\right)_0 \right]\\
& + \left( {\partial \bar{W}_x \over \partial \bar{y} } \sin\Psi_0 + {\partial \bar{W}_y \over \partial \bar{y} } \cos\Psi_0 \right)
    \left[ \left( \bar{W}_y \right)_0 + \bar{V}_0 \left( {\partial \bar{W}_y \over \partial \Delta \bar{V}_c}\right)_0 \right]\\
& + {\partial \bar{W}_x \over \partial \bar{t} } \sin \Psi_0 + {\partial \bar{W}_y \over \partial \bar{t} } \cos \Psi_0 \\
\end{split}
\end{equation}
where the partial derivatives of wind and position changes with respect to the speed increment can be obtained from Eqs. (\ref{DxDyComplete}) and (\ref{B1B2}) as
\begin{equation}
\begin{split}
\left( {\partial \bar{W}_x \over \partial \Delta \bar{V}_c}\right)_0 & = {\partial \bar{W}_x \over \partial \bar{x} } \left( {\partial \Delta \bar{x} \over \partial \Delta \bar{V}_c}\right)_0 + {\partial \bar{W}_x \over \partial \bar{y} } \left( {\partial \Delta \bar{y} \over \partial \Delta \bar{V}_c}\right)_0\\
\left( {\partial \bar{W}_y \over \partial \Delta \bar{V}_c}\right)_0 & = {\partial \bar{W}_y \over \partial \bar{x} } \left( {\partial \Delta \bar{x} \over \partial \Delta \bar{V}_c}\right)_0 + {\partial \bar{W}_y \over \partial \bar{y} } \left( {\partial \Delta \bar{y} \over \partial \Delta \bar{V}_c}\right)_0\\
\end{split}
\end{equation}
and
\begin{equation}
\begin{split}
\left( {\partial \Delta \bar{x} \over \partial \Delta \bar{V}_c}\right)_0 & = {1 \over Q_{2D}} \left[ \left({ 2 \over \Delta \bar{t}}-{\partial{\bar{W}_{y}} \over \partial{\bar{y}}} \right) \sin\Psi_0 +  \left( {\partial{\bar{W}_{x}} \over \partial{\bar{y}}} \right) \cos\Psi_0 \right] \\
\left( {\partial \Delta \bar{y} \over \partial \Delta \bar{V}_c}\right)_0 & =  {1 \over Q_{2D}} \left[ \left({\partial{\bar{W}_{y}} \over \partial{\bar{x}}} \right) \sin\Psi_0 +  \left( { 2 \over \Delta \bar{t}}-{\partial{\bar{W}_{x}} \over \partial{\bar{x}}}  \right) \cos\Psi_0 \right] \\
\end{split}
\end{equation}

\subsection{Heading Strategy}

Similarly, the incremental heading change can be obtained from Eq. (\ref{PowerApprox}) as
\begin{equation}
\Delta \Psi_c = \begin{cases} 0 & \text{$|I_{\Delta \Psi}| < \epsilon$} \\
\eta \max\{-\Delta \Psi_{c,\max}, \Psi_{\min} - \Psi_0 \} &\text{$I_{\Delta \Psi} \geq \epsilon > 0$} \\
\eta \min \{ \Delta \Psi_{c,\max},  \Psi_c \leq \Psi_{\max} - \Psi_0 \} &\text{$I_{\Delta \Psi} \leq -\epsilon $}
\end{cases}
\end{equation}
where $\bar{I}_{\Delta \Psi_c}$ could be expressed in compact form as
\begin{equation}\label{ItDeltatPsi}
\bar{I}_{\Delta \Psi_c} = \bar{V}_0 \left[ {\partial \bar{W}_V(t_0+\Delta t) \over \partial (\Delta \Psi_c)} \right]_0
\end{equation}
and

\begin{figure*}[!t]
\normalsize
\begin{equation}\label{ItDeltatPsi2}
\begin{split}
\left[ {\partial \bar{W}_V(t_0+\Delta) \over \partial (\Delta \Psi_c)} \right]_0 &
= \left( {\partial \bar{W}_x \over \partial \bar{y}} + {\partial \bar{W}_y \over \partial \bar{x}} \right) \bar{V}_0 \cos 2\Psi_0 + \left( {\partial \bar{W}_x \over \partial \bar{x}} - {\partial \bar{W}_y \over \partial \bar{y}} \right) \bar{V}_0 \sin 2\Psi_0 \\
& + \left[ {\partial \bar{W}_x \over \partial \bar{x}} \left( {\partial \Delta \bar{x} \over \partial \Delta \Psi_c}\right)_0 + {\partial \bar{W}_x \over \partial \bar{y}} \left( {\partial \Delta \bar{y} \over \partial \Delta \Psi_c}\right)_0 \right]
		\left( {\partial \bar{W}_x \over \partial \bar{x} } \sin\Psi_0 + {\partial \bar{W}_y \over \partial \bar{x} } \cos\Psi_0 \right) \\
& + \left[ {\partial \bar{W}_y \over \partial \bar{x}} \left( {\partial \Delta \bar{x} \over \partial \Delta \Psi_c}\right)_0 + {\partial \bar{W}_y \over \partial \bar{y}} \left( {\partial \Delta \bar{y} \over \partial \Delta \Psi_c}\right)_0 \right]
		\left( {\partial \bar{W}_x \over \partial \bar{y} } \sin\Psi_0 + {\partial \bar{W}_y \over \partial \bar{y} } \cos\Psi_0 \right) \\
& + \left( \bar{W}_{x} \right)_0 \left( {\partial \bar{W}_x \over \partial \bar{x} } \cos\Psi_0 - {\partial \bar{W}_y \over \partial \bar{x} } \sin\Psi_0 \right)  + \left(\bar{W}_{y}\right)_0 \left( {\partial \bar{W}_x \over \partial \bar{y} } \cos\Psi_0 - {\partial \bar{W}_y \over \partial \bar{y} } \sin\Psi_0 \right) \\
& +  \left( {\partial \bar{W}_x \over \partial \bar{t} } \cos\Psi_0 \right) - \left( {\partial \bar{W}_y \over \partial \bar{t} } \sin\Psi_0 \right) \\
\end{split}
\end{equation}
\hrulefill
\vspace*{4pt}
\end{figure*}

Again, the partial derivatives of position changes with respect to the speed increment can be obtained from Eqs. (\ref{DxDyComplete}) and (\ref{B1B2}) as
\begin{equation}
\begin{split}
\left( {\partial \Delta \bar{x} \over \partial \Delta \Psi_c}\right)_0 & = {\bar{V}_0 \over Q_{2D}} \left[ \left({ 2 \over \Delta \bar{t}}-{\partial{\bar{W}_{y}} \over \partial{\bar{y}}} \right) \cos\Psi_0 -  \left( {\partial{\bar{W}_{x}} \over \partial{\bar{y}}} \right) \sin\Psi_0 \right] \\
\left( {\partial \Delta \bar{y} \over \partial \Delta \Psi_c}\right)_0 & =  {\bar{V}_0 \over Q_{2D}} \left[ \left({\partial{\bar{W}_{y}} \over \partial{\bar{x}}} \right) \cos\Psi_0 - \left( { 2 \over \Delta \bar{t}}-{\partial{\bar{W}_{x}} \over \partial{\bar{x}}}  \right) \sin\Psi_0 \right] \\
\end{split}
\end{equation}

\section{Simulation Evaluation}


In the current paper,  the {\bf Aircraft Dynamics} is modeled by the previously described point-mass equations. The UAV tracking logic, contained within the {\bf Trajectory Tracking} concept, is based on the method of feedback linearization as described below. In general, the {\bf Wind Estimation} process represents sensors and algorithms for deriving estimates of the current wind states. The current paper seeks to focus on the development of algorithms for utilizing wind energies. It is therefore assumed that accurate wind estimates can be made and are available. Future work shall consider effects of errors associated with wind measurements and estimations.

%

\subsection{Models of Closed-Loop Tracking}

It is assumed that once optimal incremental adjustments ($\Delta \bar{V}$, $\Delta \Psi$) are derived, the UAV would track these commands in their flights. Actual onboard trajectory control logics can be very complicated and can also vary from vehicle to vehicle. In this paper, the method of feedback linearization is used to develop models of actual onboard trajectory tracking logics. The point-mass dynamic model has three control variables: $\bar{T}$ (or $\bar{P}$), $C_L$ and $\mu$. Therefore, we need to develop three closed-loop trajectory control models.

Use of the feedback linearization method starts with the specification of desired closed-loop dynamics. Specifically, a desired closed-loop airspeed tracking using thrust can be specified as a first-order system.
\begin{equation}
\dot{V} + K_V (V-V_c) = 0
\end{equation}
Using normalized variables, the closed-loop thrust law can be determined from Eq. (\ref{VdotNormalized}) as
\begin{equation}
\bar{T} = -K_V t_n (\bar{V} - \bar{V}_c)  + \bar{\rho} \bar{V}^2 \left(C_{D_0} + K C_{L}^2 \right) + \sin\gamma + \bar{W}_V'
\end{equation}
where $t_n = V_n/g$ is the normalization time. Similarly for the heading control using bank angle, we have
\begin{equation}
\dot{\Psi} + K_{\Psi}(\Psi-\Psi_c) = 0 \\
\end{equation}
which leads to
\begin{equation}\label{sinMU}
\bar{\rho}\bar{V}^2C_L \sin\mu = \bar{W}_{\Psi}' -\bar{V} \cos\gamma K_{\Psi} t_n (\Psi - \Psi_c)
\end{equation}
Finally, tracking a commanded flight path angle using lift coefficient can be achieved with
\begin{equation}
\dot{\gamma} + K_{\gamma}(\gamma-\gamma_c) = 0
\end{equation}
which results in
\begin{equation}\label{cosMU}
\bar{\rho} \bar{V}^2 C_L \cos\mu = \cos\gamma - \bar{W}_{\gamma}' -\bar{V} K_{\gamma} t_n (\gamma - \gamma_c)
\end{equation}
Combining Eqs. (\ref{sinMU}) and (\ref{cosMU}), we obtain
\begin{equation}
\tan \mu = {\bar{W}_{\Psi}' -\bar{V} cos(\gamma) K_{\Psi} t_n (\Psi - \Psi_c) \over \cos\gamma - \bar{W}_\gamma' -\bar{V} K_{\gamma} t_n (\gamma - \gamma_c)}
\end{equation}
and corresponding $C_L$ is given by Eq.(\ref{CL_control}).

\begin{figure*}[!t]
\normalsize
\begin{equation}
\begin{split}\label{CL_control}
& C_L = {\sqrt{\left[ \bar{W}_{\Psi}' -\bar{V} \cos\gamma K_{\Psi} t_n (\Psi - \Psi_c) \right]^2 + \left[ \cos\gamma - \bar{W}_{\gamma}' -\bar{V} K_{\gamma} t_n (\gamma - \gamma_c) \right]^2} \over \bar{\rho}\bar{V}^2}
\end{split}
\end{equation}
\hrulefill
\vspace*{4pt}
\end{figure*}

In the above, the feedback gains $(K_V, K_\Psi, K_\gamma)$ can be selected to reflect typical closed-loop UAV control characteristics. In this paper, it is assumed that $K_V = 0.5$, $K_\Psi = 0.5$, and $K_\gamma = 0.5$, all in sec$^{-1}$.

\subsection{Guidance Algorithm Parameters}

Performances of the proposed guidance strategies strongly depend on the following four parameters
\begin{equation}
\Delta \bar{V}_{c,\max},~~\Delta \bar{V}_{c,\min},~~\Delta \Psi_{c,\max},~~\Delta \Psi_{c,\min}
\end{equation}
Ranges of their appropriate values are now estimated.

From Eq. (\ref{VdotNormalized})
\begin{equation}\label{Vbarapprox}
\bar{V}' \approx {\Delta \bar{V} \over \Delta \bar{t}}  = {\bar{P} \over \bar{V}}- \bar{\rho} \bar{V}^2 (C_{D_0} + K C_L^2)
\end{equation}
For a steady state level flight, Eq. (\ref{CLinlevel}) suggests that $C_L \approx {1  / \bar{\rho}\bar{V}^2}$. Then, we have
\begin{equation}
\begin{split}
\Delta \bar{V}_{c,\max} \leq &~ \Delta \bar{t} \left( {\bar{P}_{\max} \over \bar{V}} - \bar{\rho} \bar{V}^2 C_{D_0} - {K \over \bar{\rho} \bar{V}^2} \right)\\
\Delta \bar{V}_{c,\min} \geq &~ \Delta \bar{t} \left( {\bar{P}_{\min} \over \bar{V}} - \bar{\rho} \bar{V}^2 C_{D_0} - {K \over \bar{\rho} \bar{V}^2} \right)\\
\end{split}
\end{equation}

Similarly for the maximum heading angle change, we have
\begin{equation}\label{Psiapprox}
\begin{split}
\Delta \Psi_{c,\max}  \leq \Delta \bar{t} \left( {1 \over \bar{V}} \sin \mu_{\max} \right) \\
\Delta \Psi_{c,\min}  \geq -\Delta \bar{t} \left( {1 \over \bar{V}} \sin \mu_{\max} \right)
\end{split}
\end{equation}

Actual values used in the guidance strategies can be smaller than the above bounds. In this paper, we select $\Delta \bar{V}_{c,\max} \leq 5$ ft/sec, $\Delta \bar{V}_{c,\min}  = - \Delta \bar{V}_{c,\max}$, $ \Delta \Psi_{c,\max} \leq 30^\circ$, and
$\Delta \Psi_{c,\min} = - \Delta \Psi_{c,\max}$.

\subsection{Evaluation Criterion}

Because the guidance strategies in this paper are introduced to save power consumptions in UAV flights, a basic performance measure is defined as the average power consumption over a specified time interval, where $t_f$ is the time period of evaluation
\begin{equation}\bar{P} = {1 \over \bar{t}_f} \int^{\bar{t}_f}_{0} \bar{T} \bar{V} d \bar{t} \end{equation}
In the following numerical results, the number of integration steps in each $[t, t+\Delta t]$ is 50, and $t_f = 50 \Delta t$.

\begin{figure}[htbp]
\centering
\includegraphics[width=4in]{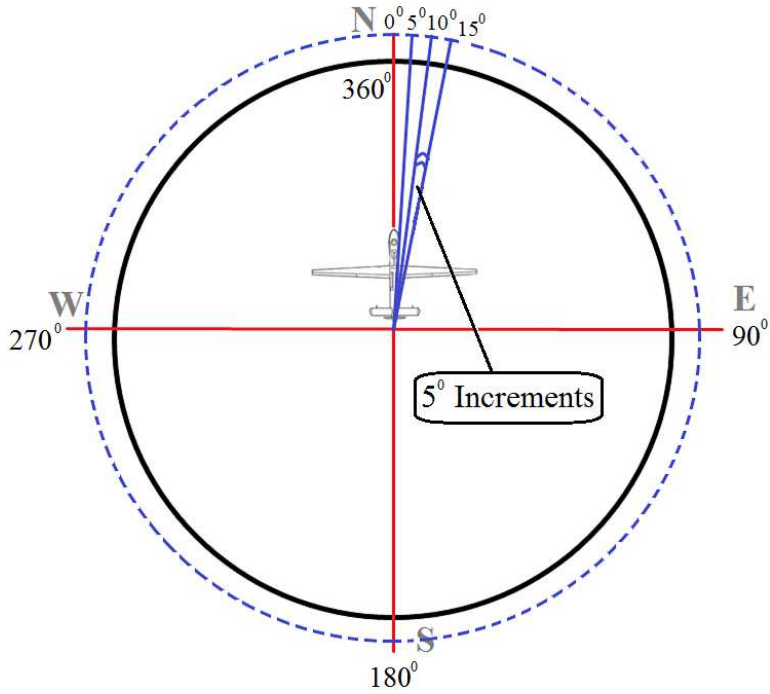}
%
%
\caption{Average power concept for different initial heading angle ($\Psi_0$) settings.}
\label{fig:AvgPsi}       
\end{figure}

Furthermore, different initial heading angles of the UAV result in different relative angles with respect to the wind field, and thus affect the resulting power saving benefit. In order to filter out these differences caused by different initial headings of UAV flights, the above basic performance measure is further averaged over a set of different initial heading angles over $[0, 360^o]$, where a generic case of $\Delta \Psi_0 = 5^{\circ}$ increment is illustrated in Fig. \ref{fig:AvgPsi}. The mean of basic average power consumptions over different initial heading conditions is defined as the measure of the performance
\begin{equation}\bar{P}_{\rm avg} = {1 \over N_\Psi} \sum^{N_\Psi}_{i=1} \bar{P}_{\Psi_0} \end{equation}
where each $i$ corresponds to a different initial heading angle, and $N_{\Psi}$ is the number of different initial heading angles used.

 \begin{figure}[htbp]
 \centering
 \includegraphics[scale=0.475]{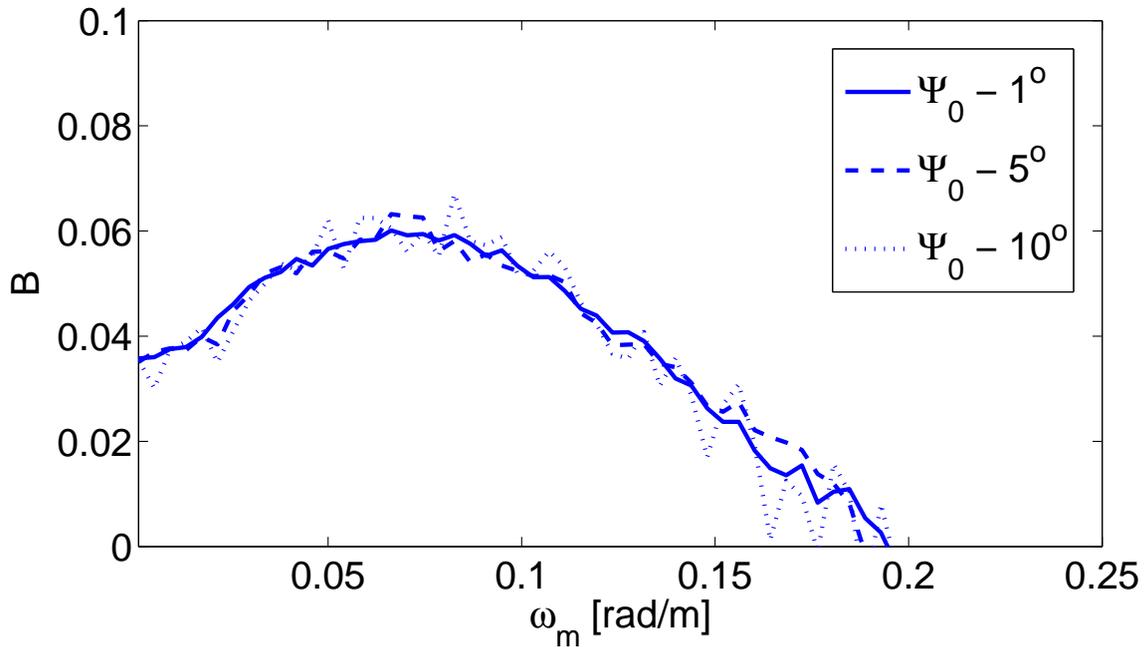}
 %
 %
 \caption{Effects of $\Psi_0$ sampling sizes for Strategy-4 ( $\Delta V$ \& $\Delta \Psi$).}
 \label{fig:Psi_0}       
 \end{figure}

In the simulation studies, effects of different initial heading angle increments are examined. Fig. \ref{fig:Psi_0} shows that $\Delta \Psi_0 = 5^o$ provides sufficiently accurate solutions as those with $\Delta \Psi_0=1^o$ at less computational loads, whereas it provides more accurate results than those with $\Delta \Psi_0=10^o$. Therefore in this paper, it is assumed that $\Delta \Psi_0=5^o$ for $\Psi_0 \in [0, 360^o]$.

To evaluate the proposed guidance strategies, four scenarios are considered. In all cases, airspeed and heading angle commands are tracked via closed-loop control logics derived above.

\begin{itemize}
\item The reference strategy that seeks to follow the reference airspeed and a constant heading angle command set at the initial heading angle. This provides the reference average power consumption $\bar{P}^0$.
\item In the second scenario, the commanded airspeed is adjusted periodically based on the current wind measurements wheres the heading angle command is set constant at the initial heading angle. The resulting average power consumption is denoted as $\bar{P}^1$.
\item In the third strategy, the heading angle is adjusted periodically based on the current wind measurements whereas the airspeed command is the same as the reference airspeed. The resulting average power consumption is denoted as $\bar{P}^2$.
\item In the fourth scenario, both airspeed and heading angle commands are adjusted periodically based on the current wind measurements. The resulting average power consumption is $\bar{P}^3$.

\end{itemize}

Finally, the following benefit criterion is introduced as a relative measure of potential fuel savings of the proposed guidance strategies over the reference strategy.
\begin{equation} B_i = {\bar{P}^0_{\rm avg} - \bar{P}^i_{\rm avg} \over \bar{P}^0_{\rm avg} } \hskip 10pt i = 1,2,3 \end{equation}

\section{Wind Field Model}

Actual wind field can be enormously complex and in general defies simple analytical models. For the convenience of studies in this paper,
wind magnitude and direction are modeled separately, and then the East and North wind components are determined from
\begin{equation}
\begin{split}
\bar{W}_x =&~ \bar{W}_m (\bar{x}, \bar{y}, t) \sin \Psi_{w} (\bar{x}, \bar{y}, t) \\
\bar{W}_y =&~ \bar{W}_m (\bar{x}, \bar{y}, t) \cos\Psi_{w} (\bar{x}, \bar{y}, t) \\
\end{split}
\end{equation}
In this paper, the following wind magnitude function is used and constant wind direction is assumed.
\begin{equation}
W_m = W_{m_0} \left[ 1 + a_x \sin(\omega_{m_x} x) + a_y \sin(\omega_{m_y} y)   \right]  \hskip 10pt \Psi_w = 90^o
\end{equation}
where the wind profile consists of a constant term plus sinusoidal components in both x and y directions. In the simulation studies below, a typical wind magnitude of 9.5 [m/sec] (or $\sim$ 21mph) is used. Furthermore, it is assumed that $a_x = a_y$ and $\omega_{m_x} = \omega_{m_y} =\omega_{m}$.

\section{Numerical Results}
UAV parameters similar to those of the ScanEagle UAV are used in the simulation studies, which has mass $m = 20$kg, reference area $S=0.55$m$^2$, parasite drag coefficient of $C_{D_0} = 0.03$, aerodynamic efficiency of $E_{\max} = 12$, maximum power available of $P_{\max} = 1,400$W, and a maximum speed of $V_{\max} = 41$ m/s. In generating numerical results, it is also assumed that
the maximum lift coefficient $C_{L_{\max}} = 1.5$, the minimum lift coefficient $C_{L_{\min}} = 0$, and the maximum bank angle range $\mu_{\max} = 40^o$. In addition, the power available is assumed to be able to vary instantaneously.

Furthermore, the minimum airspeed constraint is selected to be closed to the stall speed, whereas the maximum speed is selected to be closed to a typical cruise speed, with allowances for transient dynamics in both cases.

\begin{figure}[htbp]
 \centering
 \includegraphics[scale=0.32]{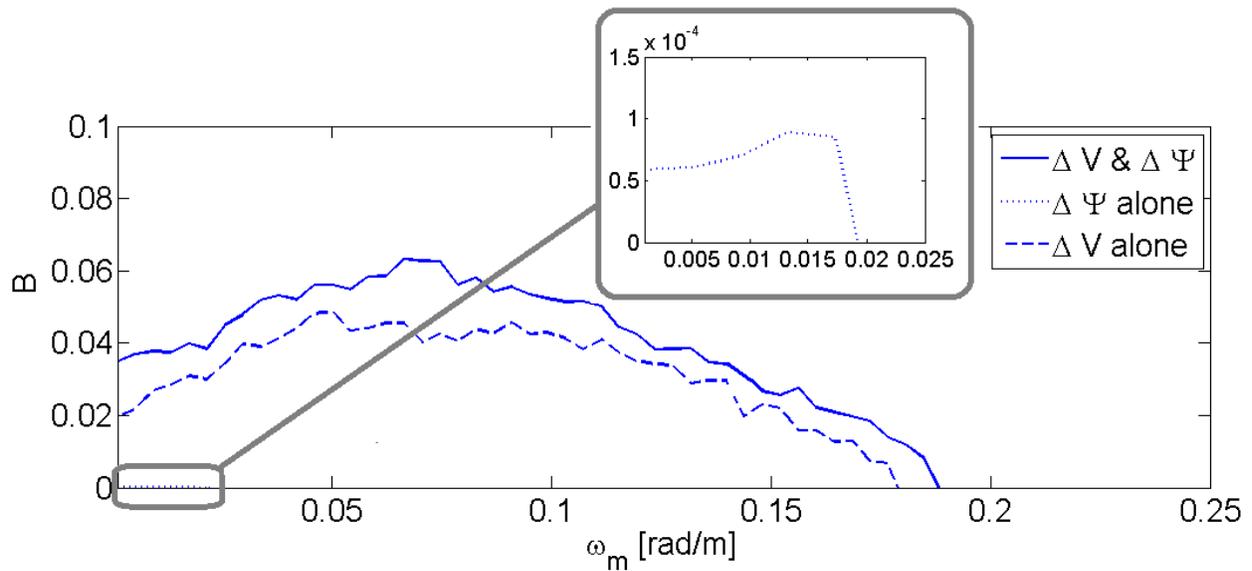}
 %
 %
 \caption{Comparison of relative benefits of the three guidance strategies, $\Delta t = 10$s.}
 \label{fig:RB1}       
\end{figure}
Fig. \ref{fig:RB1} compares the relative benefits of the three strategies in which periodic adjustments are made in the commands of airspeed alone (\emph{dashed line}), heading alone (\emph{dotted line}), and both airspeed and heading (\emph{solid line}), respectively, over the reference strategy.

Even in a constant wind field, the strategies of periodically varying airspeed and/or heading angle perform better than the constant airspeed reference trajectory. These benefits initially increase as the spatial frequency of the wind field gradually increases, but reach peaks at certain frequencies and then start to decrease beyond these frequencies. As the spatial frequency exceeds a certain limit, the proposed strategies of varying airspeed and/or heading angles derived in this paper are worse than the constant airspeed reference strategy. This is caused by the fact that as the spatial frequency exceeds a certain limit, the projected power consumption from which the airspeed and heading angle adjustments are derived starts to deviate significantly from the actual power consumption.

The most benefit seems to come from varying the airspeed. Varying both airspeed and heading angle improves the benefit further. On the other hand, varying heading angle alone produces little benefit. For numerical examples of this paper, the peak benefit by varying both airspeed and heading angle is around $5 \sim 10$\% savings of power consumptions over the reference strategy.




The choice of the update interval $\Delta t$ can (and will) directly affect the performances of the guidance strategies, which is another topic of current research. But for this study, it has been considered that the update interval is $\Delta t = 10sec$.



In general, performances of the proposed guidance strategies depend on wind field models, vehicle performance characteristics, and parameters of the proposed strategies in a complex, nonlinear way. In addition, errors in wind estimates can potentially degrade their performances. Nonetheless, results of the current paper demonstrate that the proposed real-time guidance strategies can produce positive improvements in terms of average power consumptions over the reference strategy in a fairly generic wind field. They are easy to implement and can be used in such as surveillance missions.

In practical flights, UAVs often need to stay within a certain geometric boundary due to mission as well as operational requirements. Effects of these flight constraints shall be considered in future works.

\section{Conclusions}

This paper presents real-time UAV guidance strategies that utilize wind energies to improve flight endurance. In these strategies, airspeed and/or heading angle commands are periodically adjusted based on the on-board {\sl insitu} measurements of local wind components and wind gradients. Specifically, the amounts of airspeed and/or heading angle adjustments are derived to minimize a projected power consumption at some future time. Numerical simulations are used to evaluate the relative benefits of these strategies in saving average power consumptions over a reference strategy in which the UAV follows the constant optimal steady level flight airspeed in zero wind. The average power consumption is defined over a specified interval and over different initial heading angles of the UAS. Models of onboard closed-loop tracking logics that follow adjusted commands are developed using the method of feedback linearization. UAV parameters similar to those of the ScanEagle are used in numerical simulations.

The proposed strategies offer improvements over the constant airspeed reference strategy in terms of average power consumptions. These benefits initially increase as the spatial frequency of the wind field gradually increases, but reach peaks at certain frequencies and then start to decrease beyond these frequencies. As the spatial frequency exceeds a certain limit, the strategies of varying airspeed and/or heading angles become ineffective. Here, varying the airspeed saves energy, but varying both airspeed and heading angle improves the benefit even further. On the other hand, varying heading angle alone produces little benefit. For a given spatial wind frequency, the relative benefits peak for a certain time update interval and would decrease for smaller or larger update intervals. For numerical examples of this paper, the peak benefit by varying both airspeed and heading angle is around $10$\% savings of power consumptions over the reference strategy.




%
%
%
%
%
%
%
%

\end{document}